\documentclass{amsart}
\usepackage{amsfonts,amssymb,amsmath}
\newtheorem{theorem}{Theorem}[section]
\newtheorem{lemma}[theorem]{Lemma}

\theoremstyle{definition}
\newtheorem{definition}[theorem]{Definition}
\newtheorem{proposition}[theorem]{Proposition}

\theoremstyle{remark}
\newtheorem{remark}[theorem]{Remark}

\numberwithin{equation}{section}

\newcommand{\abs}[1]{\lvert#1\rvert}

\newcounter{step}

\newtheoremstyle{algorithm}{}{}{\upshape\sffamily\mdseries%
                       \setcounter{step}{0}
                       \newcommand{\step}{\par%
                           \ifcase\arabic{step}\else\vspace{-\parskip}\fi%
                           \refstepcounter{step}%
                           \indent\hangindent\parindent%
                           \llap{\arabic{step}.\enspace}%
                           \ignorespaces}}%
                           {}{\bfseries}{.}{\newline}{}

\theoremstyle{algorithm}
\newtheorem{algorithm}[theorem]{Algorithm}
\DeclareSymbolFontAlphabet{\mathbb}{AMSb}
\DeclareMathOperator{\Ker}{Ker}
\DeclareMathOperator{\End}{End}

\title[Supersingular elliptic curves and maximal quaternionic orders]{On the correspondece between supersingular elliptic curves
and maximal quaternionic orders}
\author{Juan Marcos Cervi\~{n}o}
\address{Mathematisches Institut, Georg-August-Universit\"{a}t G\"{o}ttingen, Germany.}
\email{cervino@uni-math.gwdg.de}
\urladdr{http://www.uni-gwdg.de/cervino}
\thanks{The author wishes to thank {Prof. Ulrich Stuhler} for his continuous support and encouragement,
{Prof. Rainer Schulze-Pillot} for kindly providing answers, references and the implementation of \cite{RSP:Genus}, and {Prof. Yuri Tschinkel} for support and 
the ``Gau\ss-Lab'', where these algorithms were implemented. The author was supported by a
scholarship of the German Academic Exchange Service (DAAD), Kennziffer A/01/17917.}

\begin{document}
\begin{abstract}
We present a deterministic and explicit algorithm for the computation of the endomorphism rings of supersingular 
elliptic curves. Given any prime characteristic $p$, the algorithm returns a list of pairs $(E,\{1,e_1(E),e_2(E),
\penalty1000e_3(E)\})$ for all supersingular elliptic curves $E$ over $\overline{\Bbb{F}_p}$, where the second coordinate
 is a base of $\End_{\overline{\Bbb{F}_{p}}}(E)\otimes\Bbb{Q}$.\\
We will give at the end a table of supersingular elliptic curves with 
their respective endomorphism rings, resembling Deuring's table \cite[257-258]{Deu1}.
\end{abstract}
\maketitle
\hspace{1.in}
\begin{section}{Introduction}
In this section we fix notation and state well known results concerning elliptic curves
over finite fields.\\

Let $k$ be a finite field with $q=p^d$ elements ($p$ prime), $E$ an elliptic curve in 
$\Bbb{P}^2(\overline{k})$, given by:
\begin{equation}\label{WE}
y^2+a_1xy+a_3y=x^3+a_2x^2+a_4+a_6;\text{ with the }a_i\text{'s}\in k;
\end{equation}

plus the only point ${\rm PO}$ of the curve laying on the line at infinity. This curve $E$ is said to be 
{\sl defined in $k$,} and we will denote this by $E/k.$ Fix ${\rm PO}$ as the neutral element
of the group structure of $E.$\par
The curve $E$ is called supersingular if it satisfies the following equivalences (see \cite{Sil} 
or \cite{Huse}):

\begin{theorem}\label{thm:SSequiv}
Let $E/k$ be an elliptic curve, and denote $[\lambda]$ the isogeny multiplication by $\lambda$, with kernel $E[\lambda]$ (not to be seen as group scheme, rather as the group of points of ``order'' 
$\lambda$). Then the following are equivalent:
\begin{enumerate}
\item $E[p^r]=\Ker([p^r])=0$ for one and hence for all $r\geq1.$
\item The map $[p]$, is purely inseparable and
$j(E)\in\Bbb{F}_{p^2}.$
\item $\End_{\overline{k}}(E)$ is a non-commutative order.
\item the function field $k(E)$ has no cyclic (separable and unramified) $p$-extensions.
\end{enumerate}
\end{theorem}

In general, the endomorphism rings of elliptic curves are:
\begin{itemize}
\item $\Bbb{Z},$
\item an order of an imaginary quadratic number field,
\item a maximal order of a quaternion algebra over $\Bbb{Q}.$
\end{itemize}

The last case can only occur when the field of definition has characteristic $p>0$
and after \eqref{thm:SSequiv} when the elliptic curve is supersingular.
This last assertion was proved in \cite{Deu2}, moreover he shows there that the quaternion algebra is
exactly $\Bbb{Q}_{\infty,p}$, the only quaternion algebra (up to isomorphism) over $\Bbb{Q}$ ramified only at
$1$ and $p$. In \cite{Deu1} he proved indeed that all maximal order types of
this algebra $\Bbb{Q}_{\infty,p}$ appear as endomorphism rings of supersingular elliptic curves over 
$\Bbb{F}_{p^2}$, and that the number of those was exactly the type number of $\Bbb{Q}_{\infty,p}$, 
resulting hence a one to one correspondence between the maximal types of this quaternionic algebra
and the supersingular elliptic curves (up to isomorphism) defined on fields of characteristic $p$.

In this article we study this correspondence and give an algorithm to
effectively determine it.\\

\end{section}

\begin{section}{Arithmetic in quaternion algebras and ternary quadratic forms}\label{sec:QA3QF}

A {\sl quaternion algebra} $\mathfrak{A}$ over $\Bbb{Q}$ is a central simple four dimensional 
$\Bbb{Q}$-algebra and we will always have in mind $\Bbb{Q}_{\infty,p}$ for some prime number $p$ 
(we refer \cite{Dickson:AaA} and \cite{Reiner:MO} for basic notions on algebras and their arithmetic). As usual, $tr$ and $nr$ denote the (reduced) trace and (reduced) norm respectively.
An {\sl order} of $\mathfrak{A}$ is a subring containing $\Bbb{Z}$, which is also 
a finitely generated free $\Bbb{Z}$-module of rank $4$. We are interested in those orders which are maximal
under inclusion, and call them {\sl maximal orders} or simply {\sl orders}.

\begin{definition}
Let $\mathcal{O}$ and $\mathcal{O}^\prime$ be two (maximal quaternionic) orders in $\mathfrak{A}$. They are of the same 
{\sl type} when $\exists\alpha\in\mathfrak{A}$ with $nr(\alpha)\neq0$ such that:
$$\mathcal{O}=\alpha^{-1}\mathcal{O}^\prime\alpha.$$
By $t(\mathfrak{A})$ we denote the number of different types in $\mathfrak{A}$ and call it the {\sl type number} 
of the algebra.
\end{definition}

Let $\mathcal{O}$ be any quaternionic order, then one can make arithmetic on that order (by studying left (right) $\mathcal{O}$-ideals),
and in particular one defines the class number of it: $h(\mathcal{O})$.
The most basic theorem on the arithmetic of quaternion algebras may be the following:

\begin{theorem}[Brandt]
The class numbers and the type numbers are finite, and $t(\mathfrak{A})\leq h(\mathcal{O})$, for any order $\mathcal{O}$ in $\mathfrak{A}$. 
Moreover, if $\mathfrak{A}=\Bbb{Q}_{\infty,p}$ one has that all class numbers are the same for all (maximal) orders, and 
hence we can speak about the {\sl class number of the algebra}, denoted by $h(\mathfrak{\Bbb{Q}_{\infty,p}})=h_p$, and denote the type number $t_p$ as well.
\end{theorem}

Indeed he gave also a formula to compute the class number, a special case of Eichler's formula (in \cite{Eich:KZ}), formulated in greater generality.
The type number for the special case of our quaternion algebras was computed in \cite{Deu:TN} and 
for Eichler orders (in particular our maximal ones) of quaternion algebras over totally real
number fields in \cite{Pizer:TN}.

\begin{subsection}{Brandt-Sohn correspondence}\label{ssec:BS}

We may refer the reader for any not explained basic definition on the theory of quadratic forms either to
\cite{Eich:QF} or to \cite{Kne:QF}. All the ternary quadratic forms we are going to work with are integral and positive definite.
In \cite{Brandt:QA}, Brandt constructed maximal orders from ternary lattices {\sl via} even Clifford algebras. His idea was then exploited by W. Sohn 
in his thesis \cite{Sohn:Diss} where he proved the following:

\begin{theorem}\cite[Satz $5.1$]{Sohn:Diss}
There exists an explicit bijection between the classes of ternary quadratic forms 
of discriminant $-p$ and the order types of the quaternion algebra $\Bbb{Q}_{\infty,p}$.
\end{theorem}

In modern language, we find the explicit formulas to compute this correspondence for example in 
\cite{Brz:DQO1}.
For any ternary quadratic form $f$, we use the Seeber's notation (see \cite{Brandt:table}):
\begin{equation}\label{dfn:Seeber}
f=\sum_{1\leq i\leq j\leq 3}a_{ij}X_iX_j,\quad a_{ij}\in\Bbb{Z},
\end{equation}
will be represented by:
\begin{equation}\label{eq:3QF}
f=\begin{pmatrix}a_{11}&a_{22}&a_{33}\\a_{23}&a_{13}&a_{12}\\&&\tau \end{pmatrix};
\end{equation}
where $\tau$\footnote{This information will not be needed before we state the algorithm, and hence may not be always written.} is the number of automorphic transformations of the quadratic form $f$, which is a divisor of $24$ (not $3$).
The numbers $N_1:=a_{11}, N_2:=a_{22}, N_3:=a_{33}$ will be called {\sl succesive minima}
of $f$.\par
Already Gau\ss\ knew how to associate a lattice to a ternary quadratic form: $f\leadsto\Lambda_f$. Brandt associated 
further to any ternary quadratic space an order in a quaternion algebra: $\Lambda\leadsto\mathcal{O}_\Lambda\subset\mathfrak{A}=\mathcal{O}_\Lambda\otimes\Bbb{Q}$. Altogether,
any ternary quadratic form gives rise to an order in a quaternion algebra over $\Bbb{Q}$: $$f\quad\leadsto\quad C_0(f),$$
where $C_0(f)$ is the well known even Clifford algebra of $f$, with basis $\{1,e_1,e_2,e_3\}$ satisfying the following equations:
\begin{equation}\label{eq:BS}
\begin{split}
e_i^2&=a_{jk}e_i-a_{jj}a_{kk},\\
e_ie_j&=a_{kk}(a_{ij}-e_k),\\
e_je_i&=a_{1k}e_1+a_{2k}e_2+a_{3k}e_3-a_{ik}a_{jk},
\end{split}
\end{equation}
with $(i,j,k)$ any even permutation of $\{1,2,3\}$.
In this way, once we have a complete list of representatives of all equivalence-classes of 
ternary quadratic forms of
discriminant $-p$, we can write down explicitly the $\Bbb{Z}$-basis of representatives of all 
different order types in $\Bbb{Q}_{\infty,p}$.\\

In order to effectively compute a representative of each equivalence class of ternary quadratic 
forms of discriminant $-p$, one follows the same idea that for binary quadratic forms, namely 
defining what's called a {reduced ternary quadratic form}; in such a way that two different reduced ternary quadratic forms should not be equivalent. 

\begin{definition}\label{def:R3QF}
A ternary quadratic form $f$ like in \eqref{eq:3QF} is {\sl reduced} iff
\begin{enumerate}
\item $f({\bf x})\geq a_{ii}\quad\forall i=1,2,3\text{ and }\forall{\bf x}\in\Bbb{Z}^3
\text{ with }\gcd(x_i,\dots,x_3)=1$;
\item $a_{12}\geq0, a_{13}\geq0\text{ and }(a_{12}=0\vee a_{13}=0)\Rightarrow a_{23}\geq0$;
\item $a_{11}=a_{22}\Rightarrow\abs{a_{23}}\leq a_{13}$;
\item $a_{22}=a_{33}\Rightarrow a_{13}\leq a_{12}$.
\end{enumerate}
\end{definition}

We refer the reader to \cite{RSP:Genus} for an explicit algorithm to compute them.
As a consequence of the Mahler-Weyl and the Minkowski's inequalities (see \cite{vdW:QF}), one obtains in our 
particular case the following so called fundamental inequality:
\begin{equation}\label{eq:MW}
a_{11}a_{22}a_{33}\leq2 d.
\end{equation} 

Then pasting this algorithm with the equations \eqref{eq:BS} above, we are able to 
explicitly compute a basis of all different order types of $\Bbb{Q}_{\infty,p}$. \\

\end{subsection}
\end{section}

\begin{section}{A result of V\'{e}lu}

In this section we recall a result of V\'{e}lu, which is used in the algorithm to 
construct separable isogenies of a given degree on an elliptic curve.\\

Since we are interested in the correspondence between supersingular elliptic curves and their endomorphism rings we can suppose that 
$p\geq11$\footnote{We actually can suppose $p\geq29$. See \eqref{rmk:Aut}.}. The reason is that the class number of the quaternion algebras ramified at $\infty$ and
at $p\in\{2,3,5,7\}$ is one (see for instance \cite{Pizer:TN}), and hence the type numbers must be also $t_2=t_3=t_5=t_7=1$ by \eqref{thm:SSequiv}.\par
We set $k=\Bbb{F}_{p^d}$, and give any elliptic curve $E/k$. Pick a (finite) rational subgroup $G\subset E(k)$. Hence we have the following isogeny between elliptic curves: $\lambda:E\rightarrow E/G$.
With this notation the following holds:

\begin{theorem}\cite{Velu}
Once the equation of $E$ and all the elements of the subgroup $G$ are known, there exists a closed formula
to compute the equation of the isogeny $\lambda$ and the defining equation of the quotient $E/G$.
\end{theorem}

\begin{remark}\label{rmk:lIso}
In particular, given any positive integer $\ell\neq p$ one can construct all separable endomorphisms of degree $\ell$,
just by finding explicitly all subgroups of $E$ of order $\ell$, and then constructing the corresponding 
isogenies to the quotient. This quotient will be again an elliptic curve, since the projection 
morphism is finite separable and unramified, therefore by the Hurwitz's formula the genus of the 
two curves must be the same. Since one obtains also the equation of the quotient, we can compute in 
particular it's $j$-invariant. Then we compare it with the $j$-invariant of the original elliptic curve $E$, and then 
the isogeny will be an endomorphism if and only if both invariants coincide.
\end{remark}

\end{section}

\begin{section}{The algorithm}

\begin{algorithm}[Deuring's Correspondence]\label{Alg}
\hfill
\begin{description}
   \item[INPUT] A prime number $p$.
   \item[OUTPUT] A list of $t_p$ elements of the form $(j_0,\{1,e_1(j_0),e_2(j_0),e_3(j_0)\})$, 
	where the first coordinate runs over all supersingular $j$-invariants of $\overline{\Bbb{F}_p}$, and the second coordinate is a $\Bbb{Z}$-generating set for the lattice corresponding
      to the maximal order $\End_{\overline{\Bbb{F}_p}}(E(j_0))$ of $\Bbb{Q}_{\infty,p}$.
\end{description}
  \step Compute all reduced ternary quadratic forms of discriminant $-p$.\\
      {\bf Put:} ${\rm 3QF}=\{f_1,\dots,f_t\}$.
   \step Compute all different types of maximal orders in $\Bbb{Q}_{\infty,p}$.\\
      {\bf Put:} ${\rm QO}=\{\mathcal{O}_1,\dots,\mathcal{O}_t\}$.
   \step Compute all the $j$-invariants of supersingular elliptic curves over $\overline{\Bbb{F}_p}$.\\
      {\bf Put:} ${\rm Lj}=\{j_1,\dots,j_{2t-h},j_{2t-h+1},\dots,j_t\}$.
   \step Compute $h^\prime$ integers, which will enable us to classify the endomorphism 
      rings.\\ {\bf Put:} $\Lambda=\{\ell_1,\dots,\ell_{h^\prime}\}$.\par
      Compute all possible pairs $\left(tr(\alpha),nr(\alpha)\right)$ with $\alpha$ in any of the 
      maximal orders of ${\rm QO}$, such that $nr(\alpha)=\ell$ for some $\ell\in\Lambda$. So in this way, we will have
      for $i=1,\dots,t$:
      $${\rm SPL}_i=\{\left(tr(\alpha),nr(\alpha)\right)\mid\alpha\in\mathcal{O}_i\text{ and }nr(\alpha)=\ell\text{ for some $\ell\in\Lambda$}\}.$$
      {\bf Put:} ${\rm SPL}=\{{\rm SPL}_1,\dots,{\rm SPL}_t\}.$
   \step Construct all the separable isogenies of degree $\ell_i$ for every $i=1,\dots,h^\prime$
      on every supersingular elliptic curve $E(j)$ with $j\in Lj$ and for every one collect the pair $(tr,deg)$.
      In this way, we obtain a similar list to the one of the previous step, but on the side of elliptic curves,
      namely: $${\rm Isog}(E(j))=\{(tr(\lambda_1),deg(\lambda_1)),\dots,(tr(\lambda_{k(i,j)}),deg(\lambda_{k(i,j)}))\}.$$\\
      {\bf Put:} ${\rm Isog}=\{{\rm Isog}(E(j_1)),\dots,{\rm Isog}(E(j_t))\}.$
   \step Find the permutation $\sigma$ of $\{1,\dots,t\}$, such that ${\rm Isog}(E(j_i))$ and ${\rm SPL}_{\sigma(i)}$
      are equal as finite subsets of $\Bbb{Z}^2$.
   \step {\bf Return:} $\{\left(E(j_1),\mathcal{O}_{\sigma{1}}\right),\dots,\left(E(j_t),\mathcal{O}_{\sigma{t}}\right)\}.$
\end{algorithm}

Let us clarify the algorithm \eqref{Alg} step by step.\\

{\bf Step 1,2.} The first two steps were already done in \eqref{sec:QA3QF}.\\

{\bf Step 3.} A beautiful section of the famous paper \cite{Deu1} is exactly the computation of equations
for the $j$-invariants of the ``elliptic function fields'' which do not posses cyclic $p$-extensions (see the equivalences in \eqref{thm:SSequiv}). After studying 
first the case $p=2$, he writes explicit formulas for an invariant (denoted $A$ by Hasse) which detects
whether an elliptic function field has cyclic $p$-extensions or not, just by checking if $A$ is different or equal to
$0$, respectively. $A$ will be just a polynomial on the $j$-invariant of the elliptic function field
with coefficients depending only on $p$. We know further that all $j$-invariants of supersingular elliptic curves
are defined in $\Bbb{F}_{p^2}$, and then the polynomial $A=A(j)$ splits over $\Bbb{F}_p$ with factors 
of degree at most $2$. So we can easily compute it's roots, obtaining in this way all supersingular $j$-invariants in characteristic $p$. Finally, we put in ${\rm Lj}$ the roots of $A$ defined over 
the prime field, plus one root of each of it's quadratic factors. Since $A$ is of degree $h$, and 
there are $2t-h$ roots in the prime field, ${\rm Lj}$ has length $t$ equal to the number of
reduced ternary quadratic forms of discriminant $-p$.\\

{\bf Step 4.} This is the crucial part of the algorithm.
We base the proof of this step on an interesting theorem on ternary quadratic forms.
Let $f,g$ be ternary positive definite quadratic forms. Then we define $f\sim g\Leftrightarrow\exists T\in GL_n(\Bbb{Z})$
such that $f(T{\bf x})=g({\bf x})$, and call them {\sl integral equivalent}. For any such form $f$, we define the representation 
numbers $r(f,\ell)$ as $$r(f,\ell):=\#\{{\bf x}\in\Bbb{Z}^3:f({\bf x})=\ell\}.$$
By the main result of \cite{Schie}; the Theta-series 
of the ternary quadratic forms determine their classes. More precisely:

\begin{theorem}\cite{Schie}\label{thm:Schie}
Given two ternary quadratic forms $f$ and $g$ of discriminant $d$, there
 exists a bound $b(f)$ such that the following holds:
$$r(f,\ell)=r(g,\ell)\quad\forall\ell\leq b(f)\quad\Longrightarrow\quad f\sim g;$$
with $b(f)=min\{-1/14N_1+18/7N_2+N_3,3/2N_1-5/6N_2+17/6N_3,13/5N_1+N_2+N_3,7/2N_3\}$ and the $N_i$'s the successive minima of $f$.
\end{theorem}

So, the representation numbers classify the classes of ternary quadratic forms. Before we use 
this theorem we prove this:

\begin{lemma}\label{lemma:0}
Let $p$ be any prime and $f$ any reduced ternary quadratic form. By the Brandt-Sohn correspondence,
$f$ has associated a maximal order $\mathcal{O}_{(f)}$ in $\Bbb{Q}_{\infty,p}$. Then the form 
$(tr^2-4nr)\mid_{\mathcal{O}_{(f)}}\cong 0\perp\tilde{q}(f)$, with $\tilde{q}(f)$ positive definite and $3$-dimensional.
\end{lemma}
\begin{proof}
Follows by straightforward computation using the formulas \eqref{eq:BS}.
We just write here the resulting matrix of $\tilde{q}(f)$: 
$$\begin{pmatrix}
a_{23}^2-4a_{22}a_{33}&2a_{12}a_{33}-a_{13}a_{23}&2a_{13}a_{22}-a_{12}a_{23}\\
2a_{12}a_{33}-a_{13}a_{23}&a_{13}^2-4a_{11}a_{33}&2a_{11}a_{23}-a_{12}a_{13}\\
2a_{13}a_{22}-a_{12}a_{23}&2a_{11}a_{23}-a_{12}a_{13}&a_{12}^2-4a_{11}a_{22}\\
\end{pmatrix}$$
\end{proof}

\begin{definition}
Let $\mathcal{O}$ be an order in $\Bbb{Q}_{\infty,p}$. Define:
$$\Gamma_{b}(\mathcal{O}):=\{(tr(\alpha),nr(\alpha))\in\Bbb{Z}^2\mid\alpha\in\mathcal{O}\text{ and }nr(\alpha)\leq b\}.$$
\end{definition}

We now study the sets of the definition above for the orders $\mathcal{O}_{f}$ associated to any ternary quadratic form $f$ as in \eqref{ssec:BS} by stating our key:

\begin{proposition}\label{prp:Gamma}
Let $\{f_1,\dots,f_t\}$ be a set of representatives of reduced ternary quadratic forms. Then the sets 
$\Gamma_{b}(\mathcal{O}_{f_i})\subset\Bbb{Z}^2$ for $i=1,\dots,t$ are all different, Namely, these subsets 
characterize uniquely the types of maximal orders in $\Bbb{Q}_{\infty,p}$ and $b=O(p)$.
\end{proposition}
\begin{proof}
Let $f_1$ and $f_2$ be any two different ternary quadratic forms of discriminant $-p$. Then we set 
$q_i:=(tr^2-4nr)\mid_{\mathcal{O}_{f_i}}; i=1,2$. From \eqref{lemma:0} we have 
$q_i\cong0\perp\tilde{q_i} (i=1,2)$, with $\tilde{q_i}$ $3$-dimensional. Applying the result 
of \eqref{thm:Schie} to $f_1$ and $f_2$, we see that these two ternary quadratic forms are equivalent iff they have the same representation numbers, hence beeing the two corresponding orders of the same type. In sum, for $f_1\not\cong f_2$ there exists a bound $b_{12}$, such that for some $\ell\leq b_{12}; r(f_1,\ell)\neq r(f_2,\ell)$. Hence, in particular, the finite sets $\Gamma_{b_{12}}(\mathcal{O}_{f_1})$ and 
$\Gamma_{b_{12}}(\mathcal{O}_{f_2})$ must be different. Setting $b:=max\{b_{ij}\mid i,j=1,\dots,t\}$, the result follows directly from \eqref{thm:Schie} and Mahler-Weyl's bound \eqref{eq:MW}, which goes like $O(p)$.
\end{proof}

After \eqref{prp:Gamma} there are, say, $h^\prime$ different norms: $\Lambda:=\{\ell_1,\dots,\ell_{h^\prime}\}$ which allow us to identify
between all different maximal order types of $\Bbb{Q}_{\infty,p}$. 
More explicitly, define:
$${\rm SPL_i}:=\{(tr(\alpha),nr(\alpha))\mid\alpha\in\mathcal{O}_{f_i}\text{ and }nr(\alpha)\in\Lambda\}.$$
Then ${\rm SPL_i}={\rm SPL_j}\Leftrightarrow i=j$, and this is what we mean by {\sl indentifying all
different maximal orders}.\\

\begin{remark}
For the concrete implementation of the construction of $\Lambda$, we actually have to ``weight''
the $\ell$'s which identify different orders. The point is that the number of different subgroups
of a given order $\ell$ of the elliptic curve depends also on the number of it's factors. We 
prefer to get prime $\ell$'s, then we would have to compute just $\ell+1$ subgroups. Suppose for 
example $\ell=6$. Then we have to compute $30$ subgroups; instead for $\ell=7$ there are just $8$.
For the table \eqref{table:DeuT}, the set $\Lambda$ of {\bf step 4} could be chosen always as a 
subset of $\{3,5,7\}$, therefore having to compute at most only $18$ subgroups of each 
supersingular elliptic curve which can be made almost instantly.
\end{remark}

{\bf Step 5.}

First set $i=1,k=1$.\par
We put $\ell:=\ell_i\in\Lambda$ and also $E:=E(j_k); j_k\in{\rm Lj}$ (output of {\bf step 3}), any elliptic curve with $j$-invariant $j_k$. 
Search for the extension $\Bbb{F}_{p^{2d_\ell}}$ (say of cardinal $q_\ell$) over $\Bbb{F}_{p^2}$ produced by the points of order $\ell$ of $E$. Consider now $G:=E(\Bbb{F}_{q_\ell})$, which will be a product of two cyclic 
groups of the same order $n$ (see \cite{Schoof:GS}), so that $n^2=1+q_\ell-tr_\ell$. We must find 
two generators $P_1, P_2$ of $G$, so called {\sl echelon} generators. 
Then with these two echelon generators, one can construct all subgroups of order $\ell$ and the respective quotients as explained in \eqref{rmk:lIso}. For each quotient which is actually an endomorphism, we may check it's trace, just by looking at the respective possible traces given by the quadratic forms, and then put it in the formula $\phi^2-\left[tr(\phi)\right]\phi+\left[\ell\right]=0$, and so we find also the traces.
Obviously, we can compute them exactly, just by doing it on the Tate module, and we know to which level we should compute, just by the Hasse bound (but we practically avoided this way, as we have just explained). After all this, we have a list of the endomorphisms of the supersingular elliptic curve $E$ with it's respective traces (the degrees are all $\ell$, clear). Append this pairs
$(tr,\ell)$ to ${\rm Isog}(E)$.\par
Before returning to the begining of the previous paragraph, if $k=t$ stop, otherwise $i=i+1$. If $i>h^\prime$ put $i=1$ and $k=k+1$. Return to the begining of the previous paragraph.\\

{\bf Step 6.}

Once we computed ${\rm Isog}$, what's left is only to search the only bijection\footnote{See \eqref{rmk:Aut}.}  from ${\rm Isom}$ to
${\rm SPL}$, which exists by Deuring's correspondence. This bijection give us the permutation
$\sigma$ we are looking for, which so far means a correspondence between the supersingular elliptic 
curves and the reduced ternary quadratic forms.\\

{\bf Step 7.} With the permutation $\sigma$ and the Brandt-Sohn correspondence explained
in \eqref{ssec:BS}, we obtain the desired output.\hfil$\qed$\\

\begin{remark}\label{rmk:Aut}
We must also mention a fact which enables us also to simplify our 
previous algorithm; namely the automorphisms of the elliptic curves on the one side, and the
automorphic transformations of the reduced ternary quadratic forms on the other. The number of these
automorphic transformations is $\tau$ of \eqref{dfn:Seeber}. As one sees in the table 
\eqref{table:DeuT}, the 
$\tau$ invariants are $1,2,4$ or $6$ depending on well known facts on the $j$-invariants of the
corresponding supersingular elliptic curves (see \cite{Sil}). The ``strange'' case $\tau=1$
occurs only when there are two different isomorphism-classes of supersingular elliptic curves 
with the same endomorphism ring (i.e. the supersingular $j$-invariant is not defined in the prime
field). See \cite[Theorem 4.5]{Wat}.\par
Therefore, in order to establish the correspondence, we can suppose $p\geq29$ since in all
previous cases either the type number was one or the different $\tau$'s 
of the ternary quadratic forms were enough to decide the correspondence.
\end{remark}

\begin{remark}\label{rmk:Compl}
About the complexity of this algorithm, we can only say, that since the bound of \eqref{thm:Schie}
has order $O(p)$, it could (at least theoretically, and we guess {\sl only}) happen, that the 
$\ell_i$'s of {\bf step 4} were all prime numbers smaller than $p$. Hence we should make
$O(\pi(p))$ extensions of $\Bbb{F}_{p^2}$, and for each find echelon generators of the rational 
points. So, without sharpening the growth of the bound in \eqref{thm:Schie} we cannot expect to
have a polynomial complexity.
\end{remark}

\end{section}

\begin{section}{Table}

As an example of our algorithm we compute the correspondence between supersingular $j$-invariants
of characteristic $p$ and reduced ternary quadratic forms of discriminant $-p$, for 
$29\leq p\leq 97$, resembling Deuring's table \cite[257-258]{Deu1}.
By the explicit Brandt-Sohn correspondence \eqref{ssec:BS} one can directly compute the basis 
of the corresponding endomorphism rings with rational coefficients. We prescind from doing this
to obtain a readable table. For every prime $p$ we write in bold on each row, the supersingular
$j$-invariants of characteristic $p$. Below each of these invariants we write the corresponding
reduced ternary quadratic form in Seeber's notation \eqref{eq:3QF}.

\begin{table}
\caption{}\label{table:DeuT}
\renewcommand\arraystretch{1}
\noindent\[\tiny
\begin{array}{|c|cccc|c|c|}
\hline
p&&0\neq j\neq1728&&&1728&0\\
\hline
&&&&&&\\
29&&&{\bf 2}&{\bf 25}&{\bf \star}&{\bf 0}\\
&&&\begin{pmatrix}1&3&3\\2&0&1\\&&2\end{pmatrix}
&\begin{pmatrix}1&2&4\\1&1&0\\&&2\end{pmatrix}
&&\begin{pmatrix}1&1&10\\0&1&1\\&&6\end{pmatrix}\\
\hline
31&&&{\bf 2}&{\bf 4}&{\bf 23}&\star\\
&&&\begin{pmatrix}1&2&4\\1&0&0\\&&2\end{pmatrix}
&\begin{pmatrix}1&2&5\\2&0&1\\&&2\end{pmatrix}
&\begin{pmatrix}1&1&8\\0&1&0\\&&4\end{pmatrix}&\\
\hline
37&&&{\bf 8}&{\bf 3\pm10\sqrt{2}}&\star&\star\\
&&&\begin{pmatrix}1&2&5\\2&2&3\\&&2\end{pmatrix}
&\begin{pmatrix}2&2&3\\0&2&1\\&&1\end{pmatrix}&&\\
\hline
41&&{\bf 3}&{\bf 28}&{\bf 32}&\star&{\bf 0}\\
&&\begin{pmatrix}1&3&4\\2&1&0\\&&2\end{pmatrix}
&\begin{pmatrix}1&2&6\\1&0&1\\&&2\end{pmatrix}
&\begin{pmatrix}1&3&4\\0&1&1\\&&2\end{pmatrix}
&&\begin{pmatrix}1&1&14\\0&1&1\\&&6\end{pmatrix}\\
\hline
43&&&{\bf 41}&{\bf 12\pm8\sqrt{2}}&{\bf 8}&\star\\
&&&\begin{pmatrix}1&3&4\\1&0&1\\&&2\end{pmatrix}
&\begin{pmatrix}2&2&3\\0&1&1\\&&1\end{pmatrix}
&\begin{pmatrix}1&1&11\\0&1&0\\&&4\end{pmatrix}&\\
\hline
47&&{\bf 9}&{\bf 10}&{\bf 44}&{\bf 36}&{\bf 0}\\
&&\begin{pmatrix}1&2&7\\0&1&1\\&&2\end{pmatrix}
&\begin{pmatrix}1&2&6\\1&0&0\\&&2\end{pmatrix}
&\begin{pmatrix}1&3&4\\1&0&0\\&&2\end{pmatrix}
&\begin{pmatrix}1&1&12\\0&1&0\\&&4\end{pmatrix}
&\begin{pmatrix}1&1&16\\0&1&1\\&&6\end{pmatrix}\\
\hline
53&&{\bf 46}&{\bf 50}&{\bf 28\pm9\sqrt{2}}&\star&{\bf 0}\\
&&\begin{pmatrix}1&3&5\\2&1&0\\&&2\end{pmatrix}
&\begin{pmatrix}1&2&7\\1&1&0\\&&2\end{pmatrix}
&\begin{pmatrix}2&3&3\\-1&1&2\\&&1\end{pmatrix}
&&\begin{pmatrix}1&1&18\\0&1&1\\&&6\end{pmatrix}\\
\hline
59&{\bf 15}&{\bf 28}&{\bf 47}&{\bf 48}&{\bf 17}&{\bf 0}\\
&\begin{pmatrix}1&3&5\\1&0&0\\&&2\end{pmatrix}
&\begin{pmatrix}1&4&5\\4&0&1\\&&2\end{pmatrix}
&\begin{pmatrix}1&2&9\\2&0&1\\&&2\end{pmatrix}
&\begin{pmatrix}1&4&4\\1&0&1\\&&2\end{pmatrix}
&\begin{pmatrix}1&1&15\\0&1&0\\&&4\end{pmatrix}
&\begin{pmatrix}1&1&20\\0&1&1\\&&6\end{pmatrix}\\
\hline
61&{\bf 9}&{\bf 41}&{\bf 50}&{\bf 42\pm4\sqrt{2}}&\star&\star\\
&\begin{pmatrix}1&2&8\\1&1&0\\&&2\end{pmatrix}
&\begin{pmatrix}1&2&9\\0&1&1\\&&2\end{pmatrix}
&\begin{pmatrix}1&3&6\\-1&1&1\\&&2\end{pmatrix}
&\begin{pmatrix}2&3&3\\2&0&1\\&&1\end{pmatrix}&&\\
\hline
67&&{\bf 66}&{\bf 45\pm30\sqrt{2}}&{\bf 63\pm32\sqrt{2}}&{\bf 53}&\star\\
&&\begin{pmatrix}1&4&5\\3&1&0\\&&2\end{pmatrix}
&\begin{pmatrix}2&3&3\\1&0&1\\&&1\end{pmatrix}
&\begin{pmatrix}2&2&5\\0&2&1\\&&1\end{pmatrix}
&\begin{pmatrix}1&1&17\\0&1&0\\&&4\end{pmatrix}&\\
\hline
71&{\bf 17}&{\bf 40}&{\bf 41}&{\bf 48}&{\bf 24}&{\bf 0}\\
&\begin{pmatrix}1&4&5\\3&0&0\\&&2 \end{pmatrix}
&\begin{pmatrix}1&3&6\\1&0&0\\&&2\end{pmatrix}
&\begin{pmatrix}1&4&5\\2&0&1\\&&2\end{pmatrix}
&\begin{pmatrix}1&2&9\\1&0&0\\&&2\end{pmatrix}
&\begin{pmatrix}1&1&18\\0&1&0\\&&4\end{pmatrix}
&\begin{pmatrix}1&1&24\\0&1&1\\&&6\end{pmatrix}\\
&&&&{\bf 66}&&\\
&&&&\begin{pmatrix}1&4&5\\0&1&1\\&&2\end{pmatrix}&&\\
\hline
73&{\bf 9}&{\bf 56}&{\bf 39\pm5\sqrt{5}}&{\bf 8\pm37\sqrt{5}}&\star&\star\\
&\begin{pmatrix}1&3&7\\2&0&1\\&&2\end{pmatrix}
&\begin{pmatrix}1&2&11\\2&0&1\\&&2\end{pmatrix}
&\begin{pmatrix}2&3&4\\-1&1&2\\&&1\end{pmatrix}
&\begin{pmatrix}2&2&5\\0&1&1\\&&1\end{pmatrix}&&\\
\hline
79&{\bf 15}&{\bf 17}&{\bf 21}&{\bf 64}&{\bf 69}&\star\\
&\begin{pmatrix}1&4&5\\1&0&0\\&&2\end{pmatrix}
&\begin{pmatrix}1&3&8\\3&0&1\\&&2\end{pmatrix}
&\begin{pmatrix}1&2&10\\0&1&0\\&&2\end{pmatrix}
&\begin{pmatrix}1&5&5\\4&0&1\\&&2\end{pmatrix}
&\begin{pmatrix}1&1&20\\0&1&0\\&&4\end{pmatrix}&\\
&&&&{\bf 72\pm38\sqrt{3}}&&\\
&&&&\begin{pmatrix}2&3&4\\-2&1&1\\&&1\end{pmatrix}&&\\
\hline
83&{\bf 17}&{\bf 28}&{\bf 50}&{\bf 67}&{\bf 68}&{\bf 0}\\
&\begin{pmatrix}1&3&8\\-1&1&1\\&&2\end{pmatrix}
&\begin{pmatrix}1&2&12\\1&0&1\\&&2\end{pmatrix}
&\begin{pmatrix}1&3&7\\1&0&0\\&&2\end{pmatrix}
&\begin{pmatrix}1&4&6\\3&1&0\\&&2\end{pmatrix}
&\begin{pmatrix}1&1&21\\0&1&0\\&&4\end{pmatrix}
&\begin{pmatrix}1&1&28\\0&1&1\\&&6\end{pmatrix}\\
&&&&{\bf 38\pm35\sqrt{2}}&&\\
&&&&\begin{pmatrix}2&3&4\\2&1&1\\&&1\end{pmatrix}&&\\
\hline
89&{\bf 6}&{\bf 7}&{\bf 13}&{\bf 52}&\star&{\bf 0}\\
&\begin{pmatrix}1&4&6\\1&0&1\\&&2\end{pmatrix}
&\begin{pmatrix}1&2&13\\0&1&1\\&&2\end{pmatrix}
&\begin{pmatrix}1&5&6\\5&0&1\\&&2\end{pmatrix}
&\begin{pmatrix}1&4&7\\4&0&1\\&&2\end{pmatrix}
&&\begin{pmatrix}1&1&30\\0&1&1\\&&6\end{pmatrix}\\
&&&{\bf 66}&{\bf 76\pm39\sqrt{3}}&&\\
&&&\begin{pmatrix}1&3&8\\2&1&0\\&&2\end{pmatrix}
&\begin{pmatrix}2&3&4\\0&1&1\\&&1\end{pmatrix}&&\\
\hline
97&{\bf 1}&{\bf 20}&{\bf 45\pm28\sqrt{5}}&{\bf 76\pm3\sqrt{5}}&\star&\star\\
&\begin{pmatrix}1&5&6\\-3&1&1\\&&2\end{pmatrix}
&\begin{pmatrix}1&2&14\\1&0&1\\&&2\end{pmatrix}
&\begin{pmatrix}2&3&5\\0&1&2\\&&1\end{pmatrix}
&\begin{pmatrix}2&3&5\\3&0&1\\&&1\end{pmatrix}&&\\
&&&&{\bf 81\pm22\sqrt{5}}&&\\
&&&&\begin{pmatrix}2&2&7\\0&2&1\\&&1\end{pmatrix}&&\\
\hline
\end{array}
\]
\end{table}
\end{section}

\bibliographystyle{amsplain}

\begin{thebibliography}{99}

\bibitem{Brandt:table}
Heinrich Brandt and Oskar Intrau, \emph{{Tabellen reduzierten positiver
  tern\"{a}rer quadratischer Formen}}, {Abh. der S{\"{a}}chsischen Akad. der
  Wiss. zu Leipzig} \textbf{45} (1958), no.~4.

\bibitem{Brandt:QA}
Heinrich Brandt, \emph{{Zur Zahlentheorie der Quaternionen}}, J.-ber. Deutscher
  Math. Vereing. \textbf{53} (1943), 23--57.

\bibitem{Brz:DQO1}
Juliusz Brzezinski, \emph{Definite quaternion orders of class number one},
  Journal de Th{\'{e}}orie des Nombres (1995), no.~7, 93--96.

\bibitem{Deu1}
Max Deuring, \emph{{Die Typen der Multiplicatorenringe elliptischer
  {Funktionenk\"{o}rper}}}, Abh. Math. Sem. Hamburg \textbf{14} (1941),
  197--272.

\bibitem{Deu2}
\bysame, \emph{{Theorie der Korrespondenzen algebraischerFunktionenk\"{o}rper
  II}}, J. reine angew. Math. \textbf{183} (1941), 25--36.

\bibitem{Deu:TN}
\bysame, \emph{{Die Anzahl der Typen von Maximalordnungen einer definiten
  Quaternionenalgebra mit\ primer Grundzahl}}, Jber. Deutsch. Math. Verein.
  \textbf{54} (1950), 24--41.

\bibitem{Dickson:AaA}
L.E. Dickson, \emph{Algebras and their arithmetics}, Chicago, 1923.

\bibitem{Eich:KZ}
Martin Eichler, \emph{{{\"{U}}ber Idealklassenzahl total definiter
  Quaternionenalgebren}}, Math. Zeitschr. \textbf{43} (1937), 102--109.

\bibitem{Eich:QF}
\bysame, \emph{{Quadratische Formen und orthogonale Gruppen}}, {Zweite Auflage}
  ed., {Die Grundlehren der math. Wiss.}, no.~63, Springer Verlag, 1973.

\bibitem{Huse}
Dale {Husem\"{o}ller}, \emph{Elliptic curves}, Graduate Texts in Mathematics,
  no. 111, Springer, 1986.

\bibitem{Kne:QF}
Martin Kneser, \emph{{Quadratische Formen}}, Vorlesungsbearbeitung.
  {G{\"{o}}ttingen, $1973/1974$}.

\bibitem{Pizer:TN}
Arnold Pizer, \emph{Type numbers of {Eichler} orders}, J. reine angew. Math.
  \textbf{264} (1973), 76--102.

\bibitem{Reiner:MO}
Irving Reiner, \emph{Maximal orders}, London Mathematical Society monographs,
  no.~28, Academic Press, 1975.

\bibitem{Schoof:GS}
Ren{\'{e}} Schoof, \emph{Nonsingular plane cubic curves over finite fields}, J.
  Comb. Theory Ser. A \textbf{46} (1987), 183--211.

\bibitem{Schie}
Alexander Schiemann, \emph{Ternary positive definite quadratic forms are
  determined by their {Theta} series}, Math. Annalen \textbf{308} (1997),
  507--517.

\bibitem{Sil}
Joseph~H. Silverman, \emph{The arithmetic of elliptic curves}, Graduate Texts
  in Mathematics, no. 106, Springer, 1985.

\bibitem{Sohn:Diss}
Friedhelm Sohn, \emph{{Beitr\"{a}ge zur Zahlentheorie der tern\"{a}ren
  quadratischen Formen und der Quaternionenalgebren}}, Ph.D. thesis,
  {Westf\"{a}lische Wilhelms-Universit\"{a}t zu M\"{u}nster}, 1957.

\bibitem{RSP:Genus}
Rainer Schulze-Pillot, \emph{An algorithm for computing genera of ternary and
  quaternary quadratic forms}, {Proceedings of the International Symposium on
  Symbolic and Algebraic Computation (ISSAC), Bonn} (1991), 134--143.

\bibitem{Velu}
{V{\'{e}}lu, Jacques}, \emph{{Isog\'{e}nies entre courbes elliptiques}}, C.R.
  Acad. Sc. Paris \textbf{273} (1971), 238--241.

\bibitem{vdW:QF}
Bartel~L. van~der Waerden and Herbert Gross, \emph{{Studien zur Theorie der
  quadratischen Formen}}, {Mathematische Reihe}, no.~34, {Birkh{\"{a}}user
  Verlag}, 1968.

\bibitem{Wat}
W.C. Waterhouse, \emph{Abelian varieties over finite fields}, Ann. scient.
  {\'{E}c}. Norm. Sup. \textbf{4} (1969), 521--560, PhD. Thesis.


\end{thebibliography}

\end{document}